\documentclass[12pt]{amsart}
\usepackage[a4paper,margin=2.5cm]{geometry}
\usepackage{amsmath,amssymb,amsthm}
\usepackage{booktabs,array}
\usepackage{tikz}
\usetikzlibrary{calc}
\usepackage{pgfplots}
\pgfplotsset{compat=1.17}
\usepackage[authoryear,round]{natbib}
\usepackage[colorlinks=true,linkcolor=blue,citecolor=blue,urlcolor=blue]{hyperref}

\theoremstyle{plain}
\newtheorem{theorem}{Theorem}[section]
\newtheorem{proposition}[theorem]{Proposition}
\newtheorem{lemma}[theorem]{Lemma}
\newtheorem{corollary}[theorem]{Corollary}

\newtheorem{question}[theorem]{Question}
\theoremstyle{definition}
\newtheorem{definition}[theorem]{Definition}
\newtheorem{example}[theorem]{Example}
\newtheorem{remark}[theorem]{Remark}
\newtheorem{problem}{Problem}

\newcommand{\Z}{\mathbb{Z}}
\newcommand{\Cay}{\operatorname{Cay}}
\newcommand{\EK}{\operatorname{EK}}
\newcommand{\ecyc}{\eta_{\mathrm{cyc}}}
\newcommand{\DE}{D_E}
\newcommand{\DN}{D_N}
\newcommand{\Kbar}[1]{\overline{K}_{#1}}

\newcommand{\DepositDOI}{\href{https://doi.org/10.5281/zenodo.24457627}{10.5281/zenodo.24457627}}

\makeatletter
\g@addto@macro\@setauthors{%
  \begin{center}\footnotesize
  \textit{Department of Electrical Engineering, Energy and Automation, ENSAI,}\\
  \textit{University of Ngaound\'er\'e, P.O.~Box 455, Ngaound\'er\'e, Cameroon}\\[4pt]
  \textit{E-mail addresses:} \href{mailto:fokamrigobert@gmail.com}{fokamrigobert@gmail.com} (R.~Fokam Souop),\\
  \href{mailto:laurent.bitjoka@univ-garoua.cm}{laurent.bitjoka@univ-garoua.cm} (L.~Bitjoka)
  \end{center}}
\makeatother

\begin{document}

\title[Minimum abelian Cayley hosts]{Minimum abelian Cayley hosts: complete multipartite graphs, double stars, and a lower bound for the tree problem of Babai and S\'os}

\author{Rigobert Fokam Souop}
\author{Laurent Bitjoka}

\thanks{Corresponding author: Rigobert Fokam Souop.}

\subjclass[2020]{Primary 05C25; Secondary 05C05, 05C60, 11B13, 11P70}
\keywords{Cayley graph, induced subgraph, abelian group, double star, complete multipartite graph, Kneser's theorem}

\begin{abstract}
For a finite graph $G$, let $\eta(G)$ be the least order of a finite abelian group having a Cayley graph that contains $G$ as an induced subgraph. Babai and S\'os showed that the largest value of $\eta$ on $n$-vertex graphs is of order $n^2$, and raised the possibility of substantially smaller hosts for trees. We determine $\eta$ for all complete multipartite graphs and for blow-ups of graphs meeting the Erd\H{o}s--Kelly bound. For double stars, cyclic constructions and a rounded form of Kneser's theorem give exact values when one centre carries at most three leaves and, by exact integer computation, whenever each centre carries at most $150$ leaves. Hence a linear bound for trees, if one exists, has constant at least $3$. A counting argument shows that almost all cubic graphs need hosts of order at least a constant times $n^{6/5}$, but gives no superlinear bound for trees. We also give a computer-assisted census, with checked proofs, of the $995$ connected graphs on at most seven vertices.
\end{abstract}

\maketitle

\section{Introduction}\label{sec:intro}

Every finite graph is an induced subgraph of some Cayley graph of every sufficiently large group \citep{Babai1978}. \citet{BabaiSos1985} measured how large the group must be. A Sidon set of size $n$ in a group makes every $n$-vertex graph an induced subgraph of one of its Cayley graphs; every group of order at least $cn^3$ contains such a set, and elementary abelian groups of square order at least $n^2$ do as well. In the other direction, a counting argument shows that a group of order less than $(\tfrac12+o(1))n^2$ misses some $n$-vertex graphs. We study the same question for a single graph, with the host restricted to abelian groups.

\begin{definition}\label{def:eta}
Let $G$ be a finite graph. An injection $f\colon V(G)\to\Gamma$ into a finite abelian group $\Gamma$ is an \emph{induced embedding} into the Cayley graph $\Cay(\Gamma,S)$ if, for all distinct $u,v\in V(G)$, the vertices $f(u)$ and $f(v)$ are adjacent in $\Cay(\Gamma,S)$ exactly when $uv\in E(G)$. We write $\eta(G)$ for the least $|\Gamma|$ over all finite abelian groups $\Gamma$ and connection sets $S$ admitting an induced embedding of $G$, and $\ecyc(G)$ for the same minimum over cyclic groups.
\end{definition}

The invariant $\eta$ is finite for every graph, because large cyclic groups contain Sidon sets of any prescribed size. It is related to the \emph{representation number} $\operatorname{rep}(G)$ of \citet{ErdosEvans1989}, the least $N$ such that $G$ is an induced subgraph of the unitary Cayley graph $\Cay(\Z_N,\Z_N^\times)$: there both the group and the connection set are fixed, so $\eta(G)\le\ecyc(G)\le\operatorname{rep}(G)$. Representation numbers of stars and of complete multipartite graphs were studied by \citet{AkhtarEvansPritikin2010,AkhtarEvansPritikin2012}; once the connection set is free, the answers become much simpler (Theorem~\ref{thm:multipartite}). For comparison, induced-universal graphs with unrestricted host structure are studied by \citet{Alon2017}; the problem here is different because the host must be an abelian Cayley graph. Other work on what does and does not occur inside Cayley graphs includes \citet{Spencer1983}, \citet{GodsilImrich1987}, \citet{PotechinTsang2020}, \citet{LehnerVerret2020} and, for minimal Cayley graphs, \citet{KnauerSotoGomez2025}.

\subsection*{The tree problem}
\citet[Theorem~7.1]{BabaiSos1985} proved that every group of order greater than $n^2$ has a Cayley graph containing a given $n$-vertex tree as an induced subgraph. In their Remark~7.2 they see no reason why this could not be improved, perhaps even to $n^{1+o(1)}$, and note that they do not have a non-linear lower bound. The following is the abelian, existential form of their question.

\begin{question}\label{q:trees}
Is there a constant $c$ such that $\eta(T)\le cn$ for every tree $T$ on $n$ vertices?
\end{question}

A positive answer would not give the bound of Babai and S\'os for every group, since it asks for one host of linear order. A negative answer would give the non-linear lower bound that they lacked: if $\eta(T)>cn$, then no abelian group of order at most $cn$, in particular no cyclic group of that order, has a Cayley graph containing $T$ as an induced subgraph. We do not settle Question~\ref{q:trees}. We locate a limitation of the counting method: it proves polynomial lower bounds, $cn^{6/5}$ for almost all cubic graphs and $cn^{4/3}$ for almost all graphs of a natural $2$-degenerate class, and does not yield a superlinear estimate for trees (Theorem~\ref{thm:counting} and Corollary~\ref{cor:regular}). We also show that any admissible constant satisfies $c\ge3$: the double star $D_{a,1}$, an edge with $a$ leaves attached at one end and one leaf at the other, has $\eta(D_{a,1})=3a+2$ (Corollary~\ref{cor:three}).

\subsection*{Main results}
Section~\ref{sec:prelim} collects the basic facts, including the four necessary conditions of \citet{ErdosKelly1963}, which give a lower bound $\EK(G)$. Our results are the following.

\begin{enumerate}
\item \emph{A join bound} (Theorem~\ref{thm:join}). If $G$ contains a clique of size $c$ completely joined to an independent set of size $q$, then $\eta(G)\ge (c+1)q$. This is not implied by the Erd\H{o}s--Kelly conditions.
\item \emph{Complete multipartite graphs} (Theorem~\ref{thm:multipartite}). $\eta(K_{a_1,\dots,a_k})=k\max_i a_i$, and every abelian group of that order is an optimal host.
\item \emph{Blow-ups} (Theorem~\ref{thm:blowupexact}). If $\eta(F)$ equals the linear Erd\H{o}s--Kelly floor of $F$, then $\eta$ of the blow-up of $F$ by independent sets of size $q$ is exactly $q\,\eta(F)$. Seventy-eight connected graphs on at most seven vertices meet this hypothesis, and each seeds an infinite family of exact values.
\item \emph{Double stars} (Theorems~\ref{thm:doublestar}, \ref{thm:dslower}, \ref{thm:rounded}, \ref{thm:smallb} and~\ref{thm:finite150}). For the double star $D_{a,b}$, obtained from an edge $xy$ by attaching $a$ leaves to $x$ and $b$ leaves to $y$, with $a\ge b$, the explicit constructions give an upper bound \(f(a,b)=\min\{3(a+b)-\gcd(a,b),5a\}\). The Kneser argument gives a stronger integral lower bound \(L_*(a,b)\), which equals the monotone upper envelope of $f$ for all $1\le b\le a\le150$ by a finite exact-integer calculation. In particular, $f$ itself is not always exact: $\eta(D_{9,5})=40<41=f(9,5)$. The formula $f$ is exact for $b\le3$ (Theorem~\ref{thm:smallb}) and whenever $b\mid a$ (Corollary~\ref{cor:dsdivides}). The family $D_{a,1}$ shows that a constant $c$ as in Question~\ref{q:trees}, if it exists, is at least three.
\item \emph{Counting} (Theorem~\ref{thm:counting} and Corollary~\ref{cor:regular}). For fixed $\Delta\ge3$ there is $c_\Delta>0$ such that almost every $\Delta$-regular graph has $\eta\ge c_\Delta n^{2\Delta/(\Delta+2)}$; this is $n^{6/5}$ for cubic graphs. With a connection set of bounded size, almost every such graph needs order $n^{\Delta/2}$ (Remark~\ref{rem:boundedS}).
\item \emph{Grids} (Proposition~\ref{prop:grid}). $\eta(P_r\,\square\,P_r)\le r(r+2)$ for $r\ge3$, with equality for $r=3,4$.
\item \emph{A census} (Section~\ref{sec:census}). We compute $\eta$ and $\ecyc$ for all $995$ connected graphs on at most seven vertices, and the largest value of $\eta$ on trees with at most eleven vertices. Every lower bound comes with a checked proof, and the labellings and the digests of all formulas and proofs are deposited (Section~\ref{sec:computation}). The smallest graph whose least host is not cyclic has five vertices (Proposition~\ref{prop:smallestgap}).
\end{enumerate}

Section~\ref{sec:problems} lists open problems. A companion paper \citep{FokamBitjoka2026companion} studies degree-based lower bounds for $\eta$; among other results it proves $\eta=16$ for the Petersen graph without a computer, a value that is certified here by computation (Table~\ref{tab:named}). The two papers share the computational deposit described in Section~\ref{sec:computation}.

\section{Preliminaries}\label{sec:prelim}

Graphs are finite and simple. For a graph $G$ we write $n$, $m$, $\Delta$ and $\delta$ for its numbers of vertices and edges and its maximum and minimum degrees, and $N(v)$ for the neighbourhood of a vertex $v$. Groups are finite and abelian, written additively. For a group $\Gamma$ and a set $S\subseteq\Gamma$ with $S=-S$ and $0\notin S$, the Cayley graph $\Cay(\Gamma,S)$ has vertex set $\Gamma$, and $g,h$ are adjacent when $g-h\in S$. It is $|S|$-regular; $S=\emptyset$ is allowed. A graph isomorphic to some $\Cay(\Gamma,S)$ with $\Gamma$ abelian is an \emph{abelian Cayley graph}; it need not be connected. We write groups in invariant-factor form $\Z_{d_1}\times\dots\times\Z_{d_r}$ with $d_1\mid\dots\mid d_r$.

For an injection $f\colon V(G)\to\Gamma$ put
\[
\DE(f)=\{f(u)-f(v): uv\in E(G)\},\qquad \DN(f)=\{f(u)-f(v): u\ne v,\ uv\notin E(G)\}.
\]

\begin{lemma}[Difference criterion]\label{lem:difference}
An injection $f\colon V(G)\to\Gamma$ is an induced embedding into $\Cay(\Gamma,S)$ for some $S$ if and only if $\DE(f)\cap\DN(f)=\emptyset$; in that case $S=\DE(f)$ works.
\end{lemma}

\begin{proof}
If $f$ is an induced embedding into $\Cay(\Gamma,S)$, then $\DE(f)\subseteq S$ and $\DN(f)\cap S=\emptyset$. Conversely, $\DE(f)$ is symmetric and does not contain $0$ because $f$ is injective, so $\Cay(\Gamma,\DE(f))$ is defined. For distinct $u,v$, the difference $f(u)-f(v)$ lies in $\DE(f)$ when $uv\in E(G)$, and lies in $\DN(f)$, hence outside $\DE(f)$, when $uv\notin E(G)$.
\end{proof}

This is the abelian case of \citet[Proposition~3.1]{BabaiSos1985}. The next two facts are immediate.

\begin{lemma}[Monotonicity]\label{lem:mono}
If $F$ is an induced subgraph of $G$, then $\eta(F)\le\eta(G)$ and $\ecyc(F)\le\ecyc(G)$.
\end{lemma}

\begin{lemma}\label{lem:order}
$\eta(G)\ge n$, with equality if and only if $G$ is an abelian Cayley graph.
\end{lemma}

\begin{example}\label{ex:paths}
$\eta(C_n)=n$, and $\eta(P_n)=n+1$ for $n\ge3$: the path $P_n$ is an induced subgraph of $C_{n+1}=\Cay(\Z_{n+1},\{\pm1\})$, and it is not regular, so it is not a Cayley graph.
\end{example}

\subsection{The Erd\H{o}s--Kelly floor}
Let $f$ embed $G$ into $\Cay(\Gamma,S)$. A vertex of maximum degree has $\Delta$ neighbours in the image, all in the $|S|$-regular host, so $|S|\ge\Delta$.

\begin{lemma}[Erd\H{o}s--Kelly conditions]\label{lem:EK}
Suppose that $f$ is an induced embedding of $G$ into $\Cay(\Gamma,S)$, where $|\Gamma|=n+t$ and $d=|S|$. Then $\Delta\le d\le n+t-1$, and
\begin{enumerate}
\item[(i)] $td\ge nd-2m$;
\item[(ii)] $t\ge d-\delta$;
\item[(iii)] $t^2-(d+1)t+nd-2m\ge0$;
\item[(iv)] $(n+t)d$ is even;
\item[(v)] if $G$ is regular and is not an abelian Cayley graph, then $d\ge\Delta+1$.
\end{enumerate}
\end{lemma}

\begin{proof}
Let $U=f(V(G))$ and $W=\Gamma\setminus U$, so $|W|=t$. Each $u\in U$ has $d-\deg_G(u)$ host neighbours in $W$, so there are $s:=nd-2m$ edges between $U$ and $W$. Each vertex of $W$ has at most $d$ neighbours, which gives (i). A vertex of degree $\delta$ has $d-\delta$ distinct neighbours in $W$, which gives (ii). The host graph induced on $W$ has $(td-s)/2$ edges, at most $\binom t2$, which is (iii). The host has $(n+t)d/2$ edges, which is (iv). For (v), suppose $G$ is $\Delta$-regular and $d=\Delta$. Then no vertex of $U$ has a host neighbour in $W$, so $U$ is a union of components of $\Cay(\Gamma,S)$. These components are cosets of $K=\langle S\rangle$, each inducing a copy of $\Cay(K,S)$. Hence $G$ is isomorphic to $r$ disjoint copies of $\Cay(K,S)$, that is, to $\Cay(K\times\Z_r,S\times\{0\})$.
\end{proof}

Conditions (i)--(iv) are those of \citet{ErdosKelly1963}, who showed that, for $d=\Delta$, they characterise the least number of vertices that must be added to $G$ to obtain a $d$-regular graph containing $G$ as an induced subgraph. Here we need only that they are necessary for every admissible $d$.

\begin{definition}\label{def:EK}
$\EK(G)$ is the least $N\ge n$ for which some integer $d$ with $\Delta\le d\le N-1$ satisfies (i)--(iv) with $t=N-n$.
\end{definition}

By Lemma~\ref{lem:EK}, $\eta(G)\ge\EK(G)$. Condition (v) can add more: for the Petersen graph it raises the floor from $10$ to $13$. Conditions (i) and (ii) are increasing in $d$, which gives a closed form.

\begin{corollary}[Linear floor]\label{cor:linear}
If $m\ge1$, then
\[
\eta(G)\ \ge\ L(G):=\max\Bigl\{\,2n-\frac{2m}{\Delta},\ n+\Delta-\delta\,\Bigr\}.
\]
\end{corollary}

Among the 995 connected graphs on at most seven vertices, $\EK(G)=\eta(G)$ holds for $169$; the bound $\max\{n,\lceil 2n-2m/\Delta\rceil\}$ alone, which is condition~(i), attains $\eta$ for only $77$.

\section{Lower bounds}\label{sec:lower}

\subsection{The join bound}

\begin{theorem}[Join bound]\label{thm:join}
Suppose $G$ contains a clique $C$ with $c\ge1$ vertices and an independent set $Q$ with $q$ vertices, disjoint from $C$, such that every vertex of $Q$ is adjacent to every vertex of $C$. Then $\eta(G)\ge(c+1)q$.
\end{theorem}

\begin{proof}
Let $f$ be an induced embedding of $G$ into $\Cay(\Gamma,S)$. Write $C=\{u_1,\dots,u_c\}$; translating, we may assume $f(u_1)=0$. Put $s_j=f(u_j)$ and $R_j=f(Q)-s_j$ for $1\le j\le c$. Every $x\in Q$ is adjacent to $u_j$, so $R_j\subseteq S$.

The sets $R_1,\dots,R_c$ are pairwise disjoint. Suppose $f(x)-s_i=f(y)-s_j$ with $i\ne j$. Then $x\ne y$, because $s_i\ne s_j$. Moreover $f(x)-f(y)=s_i-s_j\in S$, since $u_iu_j$ is an edge; this makes $x$ and $y$ adjacent, which is impossible in the independent set $Q$. Hence $|S|\ge cq$.

The set $f(Q)-f(Q)$ contains $0$, has at least $q$ elements, and is disjoint from $S$ because $Q$ is independent. Therefore $|\Gamma|\ge|S|+|f(Q)-f(Q)|\ge cq+q$.
\end{proof}

For $c=1$, Theorem~\ref{thm:join} gives $\eta(G)\ge 2\max_{v}\alpha(G[N(v)])$, which is the abelian form of the correspondence between induced stars and sum-free sets \citep[Proposition~7.4]{BabaiSos1985}; sum-free subsets of finite abelian groups have a large literature of their own \citep{GreenRuzsa2005}. The join bound is not implied by Lemma~\ref{lem:EK}. For $K_{1,1,q}$ with $q\ge2$ it gives $3q$, which is the exact value (Theorem~\ref{thm:multipartite}), while Corollary~\ref{cor:linear} gives only $2q+1$. In groups of prime order $p$, the Cauchy--Davenport theorem \citep[Chapter~5]{TaoVu2006} gives $|f(Q)-f(Q)|\ge 2q-1$, and the same argument yields $p\ge(c+2)q-1$.

\subsection{Counting}

\begin{theorem}\label{thm:counting}
Let $\beta>1$, and for each $n$ let $\mathcal G_n$ be a set of labelled graphs on $\{1,\dots,n\}$ with $\log_2|\mathcal G_n|\ge\beta n\log_2 n-O(n)$.
\begin{enumerate}
\item[(a)] For every $\varepsilon>0$, the proportion of $G\in\mathcal G_n$ with $\eta(G)<(\beta-1-\varepsilon)\,n\log_2 n$ tends to $0$.
\item[(b)] If, moreover, every graph in $\mathcal G_n$ has at most $\mu n$ edges, for a constant $\mu>0$, then there is $c>0$ such that the proportion of $G\in\mathcal G_n$ with $\eta(G)<c\,n^{(\beta+\mu)/(1+\mu)}$ tends to $0$.
\end{enumerate}
\end{theorem}

Since $\beta>1$, the exponent $(\beta+\mu)/(1+\mu)$ in~(b) exceeds~$1$, so for large $n$ part~(b) is stronger than part~(a) whenever it applies; part~(a) is needed only for classes whose graphs do not have $O(n)$ edges.

\begin{proof}
A labelled graph on $\{1,\dots,n\}$ with $\eta(G)\le N$ is determined by an abelian group $\Gamma$ with $|\Gamma|\le N$, a set $S\subseteq\Gamma$ and an injection $\{1,\dots,n\}\to\Gamma$. The number of abelian groups of order $M=\prod p_i^{e_i}$ is $\prod_i p(e_i)\le\prod_i 2^{e_i}\le M$, where $p(e)$ is the number of partitions of $e$. So at most $N^2\,2^N N^n=2^{N+(n+2)\log_2N}$ labelled graphs satisfy $\eta(G)\le N$. Take $N=(\beta-1-\varepsilon)n\log_2 n$. Then $\log_2 N=\log_2 n+O(\log\log n)$, and the count is at most $2^{(\beta-\varepsilon)n\log_2 n+O(n\log\log n)}$, which is $o(|\mathcal G_n|)$. This proves (a).

For (b), the count of connection sets can be sharpened. By Lemma~\ref{lem:difference} we may take $S=\DE(f)$, which is a union of at most $m\le\mu n$ classes $\{g,-g\}$; a group of order $M$ has at most $M$ such classes. With $k=\mu n$, the number of choices of $S$ is at most $\sum_{i\le k}\binom Mi\le(eN/k)^k$ for $N\ge k$; for $M\ge k$ this is the standard estimate, and for $M<k$ the sum is $2^M<2^k\le(eN/k)^k$. So at most $N^{n+2}(eN/k)^k$ labelled graphs in $\mathcal G_n$ satisfy $\eta(G)\le N$. Put $\alpha=(\beta+\mu)/(1+\mu)$ and $N=cn^\alpha$; since $\alpha>1$, the condition $N\ge k=\mu n$ holds for all sufficiently large $n$. Then
\begin{align*}
\log_2\bigl(N^{n+2}(eN/k)^k\bigr)&=\bigl(\alpha+\mu(\alpha-1)\bigr)n\log_2 n+(1+\mu)n\log_2c+O(n)\\
&=\beta n\log_2 n+(1+\mu)n\log_2 c+O(n),
\end{align*}
where the $O(n)$ term depends on $\mu$ but not on $c$. Choosing $c$ small enough makes this at most $\log_2|\mathcal G_n|-n$, so the proportion is at most $2^{-n}$.
\end{proof}

\begin{corollary}\label{cor:regular}
\begin{enumerate}
\item[(a)] For fixed $\Delta\ge3$ there is $c_\Delta>0$ such that almost every labelled $\Delta$-regular graph on $n$ vertices ($\Delta n$ even) satisfies $\eta(G)\ge c_\Delta n^{2\Delta/(\Delta+2)}$. For cubic graphs the exponent is $6/5$, and for $4$-regular graphs $4/3$.
\item[(b)] For fixed $k\ge2$, consider the labelled graphs on $\{1,\dots,n\}$ in which each vertex $i>k$ has exactly $k$ neighbours in $\{1,\dots,i-1\}$ and there are no other edges. These graphs are $k$-degenerate, and there is $c_k>0$ such that almost all of them satisfy $\eta(G)\ge c_kn^{2k/(k+1)}$.
\end{enumerate}
\end{corollary}

\begin{proof}
By \citet{BenderCanfield1978}, the number of labelled $\Delta$-regular graphs is asymptotic to $e^{(1-\Delta^2)/4}(\Delta n)!/\bigl((\Delta n/2)!\,2^{\Delta n/2}(\Delta!)^n\bigr)$, whose base-$2$ logarithm is $(\Delta/2)n\log_2 n-O(n)$. The graphs in (b) are distinct for distinct choices, and there are $\prod_{i=k+1}^{n}\binom{i-1}{k}=2^{kn\log_2 n-O(n)}$ of them. A $\Delta$-regular graph has $\Delta n/2$ edges and a graph in (b) has fewer than $kn$. Apply Theorem~\ref{thm:counting}(b) with $\beta=\mu=\Delta/2$ and with $\beta=\mu=k$.
\end{proof}

\begin{remark}\label{rem:boundedS}
The factor counting connection sets is what keeps these exponents below $2$. If the connection set has bounded size $s$, there are at most $(N+1)^s$ choices, and the same count gives $\eta\ge c\,n^{\beta}$ for almost all graphs in $\mathcal G_n$; for $\Delta$-regular graphs this is $cn^{\Delta/2}$, the order of the smallest induced-universal graph for maximum degree $\Delta$ \citep{AlonNenadov2019}. Since $\eta(G)\le4n^2$ for every graph (Section~\ref{sec:constructions}), almost every $\Delta$-regular graph with $\Delta\ge5$ has no host of order $O(n^2)$ with a connection set of bounded size. By contrast, the twisted grids of Proposition~\ref{prop:grid} embed $P_r\,\square\,P_r$, which has $n=r^2$ vertices, in a Cayley graph of order $r(r+2)$, linear in $n$, with a connection set of four elements.
\end{remark}

There are $n^{n-2}$ labelled trees, so $\beta=1$ for trees, and with $\mu=1$ the exponent in Theorem~\ref{thm:counting}(b) is $1$: the theorem is silent on Question~\ref{q:trees}, while cubic graphs and $2$-degenerate graphs already need hosts of polynomially superlinear order.

\begin{remark}[Entropy threshold and trees]\label{rem:entropy}
Theorem~\ref{thm:counting} depends on the number of \emph{labelled} graphs in the chosen class, not on the presence or absence of cycles. Its superlinear conclusion requires a counting exponent $\beta>1$ and sufficiently many edges relative to that entropy. For example, every proper minor-closed family has at most $n!c^n$ labelled graphs on $n$ vertices for a constant $c$ \citep{NorineSeymourThomasWollan2006}; the argument of Theorem~\ref{thm:counting} therefore does not automatically imply a superlinear host-order bound for such families. For trees, Cayley's formula gives only $n^{n-2}$ labelled objects, and the hypothesis $\beta>1$ fails. This says nothing about whether trees admit uniformly linear abelian hosts.
\end{remark}

\section{Constructions}\label{sec:constructions}

For a graph $F$ and $q\ge1$, the \emph{blow-up} $F[\Kbar q]$ replaces each vertex $v$ by $q$ pairwise non-adjacent copies $(v,1),\dots,(v,q)$, with $(u,i)$ adjacent to $(v,j)$ exactly when $uv\in E(F)$. The blow-up $F[K_q]$ is defined in the same way, except that the copies of each vertex form a clique.

\begin{lemma}[Blow-up]\label{lem:blowup}
Let $f$ embed $F$ into $\Cay(\Gamma,S)$, and let $Q$ be an abelian group with $|Q|\ge q$. Fix an injection $\iota\colon\{1,\dots,q\}\to Q$. Then $(v,i)\mapsto(f(v),\iota(i))$ embeds $F[\Kbar q]$ into $\Cay(\Gamma\times Q,\,S\times Q)$, and $F[K_q]$ into $\Cay(\Gamma\times Q,\,(S\times Q)\cup(\{0\}\times(Q\setminus\{0\})))$. Consequently $\eta(F[\Kbar q])\le q\,\eta(F)$ and $\eta(F[K_q])\le q\,\eta(F)$.
\end{lemma}

\begin{proof}
The image of a pair is $(f(u)-f(v),\iota(i)-\iota(j))$. If $u\ne v$, the first coordinate is non-zero and lies in $S$ exactly when $uv\in E(F)$. If $u=v$ and $i\ne j$, the pair is $(0,\iota(i)-\iota(j))$ with non-zero second coordinate: it lies outside $S\times Q$ and inside $\{0\}\times(Q\setminus\{0\})$.
\end{proof}

Two non-adjacent vertices with the same neighbourhood are \emph{open twins}; this is an equivalence relation. Let $G^\circ$ be the induced subgraph obtained by keeping one vertex from each class, and let $\tau(G)$ be the size of a largest class.

\begin{corollary}[Twin quotient]\label{cor:twins}
$\eta(G^\circ)\le\eta(G)\le\tau(G)\,\eta(G^\circ)$. The classes, and hence $G^\circ$, can be found in $O(n+m)$ expected time by hashing neighbourhoods, and an embedding of $G^\circ$ lifts to an embedding of $G$ in $O(n^2)$ group operations.
\end{corollary}

\begin{proof}
$G^\circ$ is an induced subgraph of $G$, and $G$ is an induced subgraph of $G^\circ[\Kbar{\tau}]$, since twins in one class have the same neighbours. Apply Lemmas~\ref{lem:mono} and~\ref{lem:blowup}.
\end{proof}

\begin{proposition}[Products]\label{prop:product}
$\eta(G\,\square\,H)\le\eta(G)\,\eta(H)$.
\end{proposition}

\begin{proof}
Let $f$ embed $G$ into $\Cay(\Gamma_1,S_1)$ and $g$ embed $H$ into $\Cay(\Gamma_2,S_2)$. Map $(u,v)\mapsto(f(u),g(v))$ into $\Cay(\Gamma_1\times\Gamma_2,\,(S_1\times\{0\})\cup(\{0\}\times S_2))$. Since $f$ and $g$ are injective, $(f(u)-f(u'),g(v)-g(v'))$ lies in the connection set exactly when either $uu'\in E(G)$ and $v=v'$, or $u=u'$ and $vv'\in E(H)$.
\end{proof}

For grids, Proposition~\ref{prop:product} and Example~\ref{ex:paths} give $\eta(P_r\,\square\,P_r)\le(r+1)^2$. A twisted cyclic host does one better.

\begin{proposition}[Twisted grid]\label{prop:grid}
For $r\ge3$, $\eta(P_r\,\square\,P_r)\le r(r+2)$, and $\eta(P_r\,\square\,P_r)\ge r^2+r$.
\end{proposition}

\begin{proof}
Let $N=r(r+2)$ and label $(i,j)\mapsto\lambda(i,j)=i+(r+1)j\in\Z_N$ for $0\le i,j\le r-1$. As integers these labels are distinct (base $r+1$) and lie in $[0,(r-1)(r+2)]$, so they are distinct modulo $N$. Take $S=\{\pm1,\pm(r+1)\}$; horizontal and vertical edges have differences $\pm1$ and $\pm(r+1)$. For any other pair, the integer difference is $d=\Delta i+(r+1)\Delta j$ with $|\Delta i|,|\Delta j|\le r-1$, so $|d|\le r^2+r-2<N-(r+1)$; hence $d\equiv s\pmod N$ with $s\in S$ forces $d=s$. If $d=\pm1$ and $\Delta j\ne0$, then $|\Delta i|\ge r$, which is impossible; so $\Delta j=0$ and $\Delta i=\pm1$, an edge. If $d=\pm(r+1)$, then $r+1$ divides $\Delta i$, so $\Delta i=0$ and $\Delta j=\pm1$, again an edge. The lower bound is the first term of Corollary~\ref{cor:linear}.
\end{proof}

By computation, $\eta(P_3\,\square\,P_3)=15$ and $\eta(P_4\,\square\,P_4)=24$, so Proposition~\ref{prop:grid} is exact for $r=3,4$.

Finally, every graph has a quadratic host: if $p\ge n$ is an odd prime, the parabola $\{(x,x^2):x\in\Z_p\}$ is a Sidon set in $\Z_p^2$, and by Lemma~\ref{lem:difference} any labelling of $V(G)$ by elements of a Sidon set is an induced embedding. Hence $\eta(G)\le p^2$, in line with \citet{BabaiSos1985}.

\section{Exactly solved families}\label{sec:exact}

\begin{theorem}\label{thm:multipartite}
For all $k\ge1$ and $a_1,\dots,a_k\ge1$,
\[
\eta(K_{a_1,\dots,a_k})=\ecyc(K_{a_1,\dots,a_k})=k\max_i a_i ,
\]
and every abelian group of order $k\max_i a_i$ is an optimal host.
\end{theorem}

\begin{proof}
Put $q=\max_i a_i$. Let $\Gamma$ be an abelian group of order $kq$ and $Q\le\Gamma$ a subgroup of order $q$; it exists because finite abelian groups have subgroups of every order dividing the group order. In $\Cay(\Gamma,\Gamma\setminus Q)$ two elements are adjacent exactly when they lie in different cosets of $Q$, so this graph is $K_{q,\dots,q}$ with $k$ parts. It contains $K_{a_1,\dots,a_k}$ as an induced subgraph. For $k=1$, the graph is edgeless of order $a_1$, and its least host order is $a_1$ using the empty connection set; thus the claim holds. Assume henceforth $k\ge2$. For the lower bound, let $P$ be a largest part, and choose one vertex in each of the other $k-1$ parts. These vertices form a clique completely joined to the independent set $P$, so Theorem~\ref{thm:join} gives $\eta\ge kq$.
\end{proof}

Theorem~\ref{thm:multipartite} contains $\eta(K_{1,q})=2q$ and $\eta(K_{a,b})=2\max(a,b)$, and its proof does not use the classification of maximal sum-free sets. For the representation number the corresponding problem is much harder \citep{AkhtarEvansPritikin2012}.

\begin{theorem}\label{thm:blowupexact}
Let $F$ have at least one edge. If $\eta(F)=L(F)$, then $\eta(F[\Kbar q])=q\,\eta(F)$ for every $q\ge1$, and $\Gamma\times Q$ is an optimal host for every optimal host $\Gamma$ of $F$ and every abelian group $Q$ of order $q$.
\end{theorem}

\begin{proof}
Under blow-up, $(n,m,\Delta,\delta)$ becomes $(qn,q^2m,q\Delta,q\delta)$, so $L(F[\Kbar q])=q\,L(F)=q\,\eta(F)$. Corollary~\ref{cor:linear} and Lemma~\ref{lem:blowup} give the two bounds.
\end{proof}

\begin{example}\label{ex:blowups}
For $r\ge3$, $L(P_r)=r+1=\eta(P_r)$, so $\eta(P_r[\Kbar q])=q(r+1)$; in particular $\eta(P_4[\Kbar q])=5q$. For the graph $K_4^+$ obtained from $K_4$ by adding a pendant vertex, $L(K_4^+)=5+4-1=8=\eta(K_4^+)$ (Proposition~\ref{prop:smallestgap}), so $\eta(K_4^+[\Kbar q])=8q$. Among the 995 connected graphs on at most seven vertices, exactly $78$ satisfy $\eta(F)=L(F)$, and $64$ of them are not Cayley graphs; each seeds an infinite family of exact values.
\end{example}

\section{Double stars and trees}\label{sec:doublestars}

The double star $D_{a,b}$ ($a,b\ge1$) consists of an edge $xy$, $a$ leaves attached to $x$ and $b$ leaves attached to $y$; it has $n=a+b+2$ vertices. Double stars are the trees of diameter three. Corollary~\ref{cor:linear} gives only $\eta(D_{a,b})\ge 2a+2b+2-2b/(a+1)$ when $a\ge b$, for example $4q+1$ for $D_{q,q}$ and $2a+4$ for $D_{a,1}$.

\begin{theorem}\label{thm:doublestar}
For all $a,b\ge1$,
\[
\eta(D_{a,b})\ \le\ \min\bigl\{\,3(a+b)-\gcd(a,b),\ 5\max(a,b)\,\bigr\}.
\]
\end{theorem}

\begin{proof}
Put $k=a+b$ and $N=3k-1$, and in $\Z_N$ let $S=\{k,k+1,\dots,2k-1\}$. Since $N-j$ runs over the same interval as $j$, we have $S=-S$, and $\Cay(\Z_N,S)$ joins two residues exactly when their circular distance is at least $k$. Label
\[
x\mapsto0,\qquad y\mapsto 2k-1,\qquad \text{leaves of }x\mapsto k,\dots,k+a-1,\qquad \text{leaves of }y\mapsto a,\dots,k-1 .
\]
These $k+2$ labels are distinct. The required edges have differences in $S$: from $x$ to its leaves $k,\dots,k+a-1$; from $x$ to $y$, $2k-1$; from $y$ to its leaves, $k,\dots,k+b-1$. Every other pair has a difference of circular size between $1$ and $k-1$: from $x$ to the leaves of $y$, $a,\dots,k-1$; from $y$ to the leaves of $x$, $b,\dots,k-1$; between two leaves, $1,\dots,k-1$. So the labelling is an induced embedding, and $\eta(D_{a,b})\le 3(a+b)-1$.

Let $g=\gcd(a,b)$. Then $D_{a,b}$ is an induced subgraph of $D_{a/g,b/g}[\Kbar g]$ (keep one copy of each centre), so Lemma~\ref{lem:blowup} gives $\eta(D_{a,b})\le g\,(3(a+b)/g-1)=3(a+b)-g$. Finally, let $M=\max(a,b)$. Then $D_{a,b}$ is an induced subgraph of $P_4[\Kbar M]$, and $P_4$ is an induced subgraph of $C_5$, so $\eta(D_{a,b})\le5M$.
\end{proof}

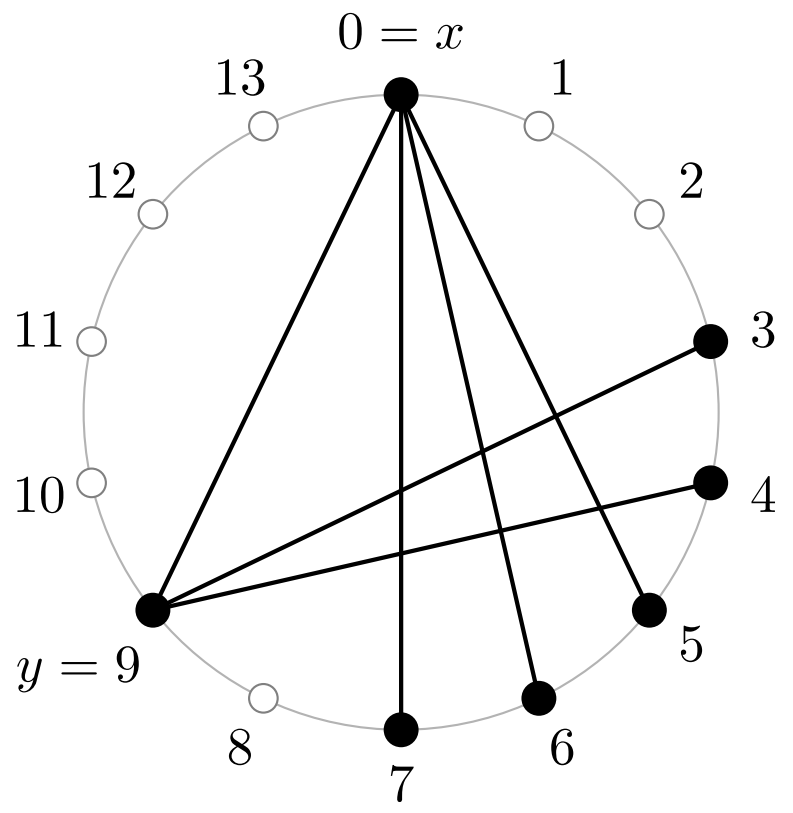
\begin{figure}[ht]
\centering
\begin{tikzpicture}[scale=2.1]
\foreach \i in {0,...,13} {
  \coordinate (p\i) at ({90-360*\i/14}:1);
}
\draw[gray!60] (0,0) circle (1);
\foreach \i/\j in {0/5,0/6,0/7,0/9,9/3,9/4} { \draw[thick] (p\i) -- (p\j); }
\foreach \i in {1,2,8,10,11,12,13} { \fill[white] (p\i) circle (0.045); \draw[gray] (p\i) circle (0.045); }
\foreach \i in {0,3,4,5,6,7,9} { \fill[black] (p\i) circle (0.055); }
\foreach \i in {1,...,8,10,11,12,13} { \node[font=\footnotesize] at ({90-360*\i/14}:1.17) {$\i$}; }
\node[font=\footnotesize] at ({90}:1.2) {$0=x$};
\node[font=\footnotesize] at ({90-360*9/14}:1.3) {$y=9$};
\end{tikzpicture}
\caption{The double star $D_{3,2}$ as an induced subgraph of $\Cay(\Z_{14},\{5,\dots,9\})$: $x\mapsto0$, $y\mapsto9$, the leaves of $x$ at $5,6,7$ and those of $y$ at $3,4$. Two residues are adjacent exactly when their circular distance is at least $5$.}
\label{fig:d32}
\end{figure}

Figure~\ref{fig:d32} shows the construction for $D_{3,2}$ in $\Z_{14}$. Every tree of diameter three on $n$ vertices therefore has $\eta\le 3n-7$.

For lower bounds we use Kneser's theorem \citep{Kneser1953}; see \citet[Chapter~5]{TaoVu2006}. For finite non-empty subsets $A,B$ of an abelian group, with $H=\{h:h+A+B=A+B\}$ the stabiliser of $A+B$, it states that $|A+B|\ge|A+H|+|B+H|-|H|$.

\begin{proposition}\label{prop:kneser}
Let $f$ embed $D_{a,b}$ into $\Cay(\Gamma,S)$, let $L$ be the set of the $a+b$ leaf labels, and let $H$ be the stabiliser of $L-L$. Then $|\Gamma|\ge3(a+b)-|H|$. In particular, every host of prime order $p$ satisfies $p\ge3(a+b)-1$, and equality is attained whenever $3(a+b)-1$ is prime.
\end{proposition}

\begin{proof}
Translate so that $f(x)=0$, and put $s=f(y)$. Let $A$ be the set of labels of the leaves of $x$ and $B'$ that of the leaves of $y$. Then $A\subseteq S$ and $B'-s\subseteq S$. These two sets are disjoint: if $a'=b'-s$ with $a'\in A$ and $b'\in B'$, then $b'-a'=s\in S$, so a leaf of $x$ and a leaf of $y$ would be adjacent. Hence $|S|\ge a+b$. The leaves form an independent set, so $L-L$ is disjoint from $S$, and $|\Gamma|\ge|S|+|L-L|$. Kneser's theorem with $A=L$ and $B=-L$ gives $|L-L|\ge2|L|-|H|$, which proves the first claim. If $|\Gamma|=p$ is prime, then $L-L\ne\Gamma$ because $S\ne\emptyset$, and the Cauchy--Davenport theorem gives $|L-L|\ge 2(a+b)-1$. Equality comes from the host in the proof of Theorem~\ref{thm:doublestar}.
\end{proof}

\begin{theorem}\label{thm:dslower}
Let $a\ge b\ge1$. Every host of $D_{a,b}$ has at least $3a+2b$ elements.
\end{theorem}

\begin{proof}
Keep the notation of the proof of Proposition~\ref{prop:kneser}, with $x$ the centre carrying $a$ leaves: $f(x)=0$, $s=f(y)\in S$, $A\subseteq S$, $B'-s\subseteq S$, $L=A\cup B'$, and $H$ is the stabiliser of $L-L$. Put $h=|H|$ and $Q=\Gamma/H$. For $X\subseteq\Gamma$ let $\bar X\subseteq Q$ be the set of cosets of $H$ that meet $X$.

Since $L-L$ is a union of cosets of $H$, each coset in $\bar L-\bar L$ is contained in $L-L$, which is disjoint from $S$. So $\bar L-\bar L$ and $\bar S$ are disjoint subsets of $Q$, and
\[
|\Gamma|=h|Q|\ \ge\ h\bigl(|\bar L-\bar L|+|\bar S|\bigr).
\]
An element of $Q$ that stabilises $\bar L-\bar L$ lifts to an element of $\Gamma$ that stabilises its preimage $L-L$, hence lies in $H$; so the stabiliser of $\bar L-\bar L$ in $Q$ is trivial, and Kneser's theorem in $Q$ gives $|\bar L-\bar L|\ge2|\bar L|-1$.

The sets $\bar A$ and $\overline{B'-s}$ are contained in $\bar S$ and are disjoint: a common coset would give $a'+H=b'-s+H$ for some $a'\in A$ and $b'\in B'$, hence $s\in b'-a'+H\subseteq L-L$, whereas $s\in S$. Therefore $|\bar S|\ge|\bar A|+|\overline{B'}|$. Finally $h|\bar L|\ge|L|=a+b$, $h|\bar A|\ge|A|=a$ and $|\overline{B'}|\ge1$, so
\[
|\Gamma|\ \ge\ h\bigl(2|\bar L|-1\bigr)+h|\bar A|+h|\overline{B'}|\ \ge\ 2(a+b)-h+a+h\ =\ 3a+2b. \qedhere
\]
\end{proof}

\begin{theorem}[Rounded Kneser obstruction]\label{thm:rounded}
For all $a,b\ge1$, write $k=a+b$ and put
\[
L_*(a,b):=\min_{1\le h\le k-1}h\left(2\left\lceil\frac{k}{h}\right\rceil-1+\left\lceil\frac a h\right\rceil+\left\lceil\frac b h\right\rceil\right).
\]
Then $\eta(D_{a,b})\ge L_*(a,b)$. Moreover $L_*(a,b)\ge3\max(a,b)+2\min(a,b)$, so Theorem~\ref{thm:rounded} contains Theorem~\ref{thm:dslower}.
\end{theorem}
\begin{proof}
Keep the notation of Theorem~\ref{thm:dslower}, with $L=A\cup B'$, $H=\operatorname{Stab}(L-L)$, $h=|H|$, $Q=\Gamma/H$, and overbars denoting images in $Q$. Since the full preimage of $\bar L-\bar L$ is $L-L$ and this set is disjoint from $S$, the subsets $\bar L-\bar L$ and $\bar S$ of $Q$ are disjoint. The difference set $\bar L-\bar L$ has trivial stabiliser: any stabilising coset would lift to an element stabilising $L-L$, and therefore belong to $H$. Kneser's theorem in $Q$ gives $|\bar L-\bar L|\ge2|\bar L|-1$. The sets $\bar A$ and $\overline{B'-s}$ lie in $\bar S$ and are disjoint, for otherwise $s\in B'-A+H\subseteq L-L$, contrary to $s\in S$. Since $\overline{B'-s}=\bar B'-(s+H)$ is a translate of $\bar B'$ in $Q$, we have $|\overline{B'-s}|=|\bar B'|$. Each coset of $H$ has exactly $h$ elements, so $h|\bar L|\ge|L|=k$, $h|\bar A|\ge|A|=a$ and $h|\bar B'|\ge|B'|=b$; as these cardinalities are integers, $|\bar L|\ge\lceil k/h\rceil$, $|\bar A|\ge\lceil a/h\rceil$ and $|\bar B'|\ge\lceil b/h\rceil$. It follows that
\begin{align*}
|\Gamma|&\ge h\bigl(2|\bar L|-1+|\bar A|+|\bar B'|\bigr)\\
&\ge h\left(2\left\lceil\frac{k}{h}\right\rceil-1+\left\lceil\frac a h\right\rceil+\left\lceil\frac b h\right\rceil\right).
\end{align*}
At $h=1$ the expression is $3k-1$. For $h\ge k$ it is at least $3h\ge3k$; hence some minimum of the displayed integer expression is attained in $1\le h\le k-1$. The actual $h$ need not be the minimiser, but the displayed inequality establishes the theorem.

For the last claim, put $M=\max(a,b)$ and $\mu=\min(a,b)$. For every $h\ge1$ we have $h\lceil k/h\rceil\ge k$, $h\lceil M/h\rceil\ge M$ and $h\lceil\mu/h\rceil\ge h$; since the expression is symmetric in $a$ and $b$, each term of the minimum is at least $2k-h+M+h=3M+2\mu$.
\end{proof}

\begin{corollary}\label{cor:dsdivides}
If $a\ge b\ge1$ and $b$ divides $a$, then $\eta(D_{a,b})=3a+2b$. In particular $\eta(D_{a,1})=3a+2$ and $\eta(D_{a,a})=5a$ for every $a\ge1$.
\end{corollary}

\begin{proof}
Here $\gcd(a,b)=b$, so Theorem~\ref{thm:doublestar} gives $\eta(D_{a,b})\le3(a+b)-b=3a+2b$, and Theorem~\ref{thm:dslower} gives the reverse inequality.
\end{proof}

\begin{theorem}[Exact values with at most three leaves at one centre]\label{thm:smallb}
For $a\ge b\ge1$ with $b\le3$,
\[
\eta(D_{a,b})=\min\{3(a+b)-\gcd(a,b),\,5a\}.
\]
\end{theorem}
\begin{proof}
The case $b=1$ is Corollary~\ref{cor:dsdivides}. If $b=2$ and $a$ is even, the universal bound of Theorem~\ref{thm:dslower} is $3a+4$, which equals the upper bound $3(a+2)-2$. Let $b=2$ and $a$ be odd ($a\ge3$). Write $F_h$ for the $h$-th expression under the minimum in Theorem~\ref{thm:rounded}. Then $F_1=3a+5$. For $h\ge2$, since $\lceil2/h\rceil=1$,
\[F_h=2h\left\lceil\frac{a+2}{h}\right\rceil+h\left\lceil\frac a h\right\rceil.
\]
If $h\nmid a$, the two terms are at least $2(a+2)$ and $a+1$, giving $F_h\ge3a+5$. If $h\mid a$, then $h$ is odd and at least three, and $h\lceil(a+2)/h\rceil=a+h$; hence $F_h=3a+2h\ge3a+6$. Thus $L_*(a,2)=3a+5$, which agrees with $3(a+2)-\gcd(a,2)$ and is at most $5a$ for $a\ge3$.

Now let $b=3$. If $3\mid a$, the lower bound $3a+6$ of Theorem~\ref{thm:dslower} agrees with the construction $3(a+3)-3$. Suppose $3\nmid a$ (so $a\ge4$). Then $F_1=3a+8$. For $h=2$ one computes $F_2=3a+10$ if $a$ is even and $F_2=3a+9$ if $a$ is odd. For $h\ge3$, the formula is
\[F_h=2h\left\lceil\frac{a+3}{h}\right\rceil+h\left\lceil\frac a h\right\rceil.
\]
If $h\mid a$, necessarily $h\ge4$ (because $3\nmid a$), and $F_h=3a+2h\ge3a+8$. Otherwise $h\lceil(a+3)/h\rceil\ge a+3$ and $h\lceil a/h\rceil\ge a+1$, so $F_h\ge3a+7$. Equality $F_h=3a+7$ would require both $h\mid(a+3)$ and $h\mid(a+1)$, hence $h\mid2$, impossible for $h\ge3$. Therefore $F_h\ge3a+8$ in all cases. We conclude $L_*(a,3)=3a+8$, equal to the upper construction $3(a+3)-1\le5a$ for $a\ge4$. These cases exhaust $b\le3$.
\end{proof}

\begin{example}[The formula $f$ is not always exact]\label{ex:counterexample}
The equality $\eta(D_{a,b})=f(a,b)$, with $f(a,b)=\min\{3(a+b)-\gcd(a,b),5a\}$, fails in general. Indeed, $f(9,5)=41$ and $f(10,5)=40$, while $D_{9,5}$ is an induced subgraph of $D_{10,5}$. Therefore $\eta(D_{9,5})\le40$. A direct certificate is the following:
\[
\Gamma=\Z_{40},\quad S=\{3,4,5\}+8\Z_{40},\quad x\mapsto0,\quad y\mapsto21,
\]
with the nine leaves of $x$ labelled $3,11,12,19,20,27,28,35,36$ and the five leaves of $y$ labelled $2,10,18,26,34$. The $16$ labels are distinct. Modulo eight the $x$-leaf labels are $3$ or $4$, the $y$-leaf labels are $2$, and $y\equiv5$; hence precisely the centre--centre and assigned centre--leaf pairs have differences in the three residue classes $3,4,5$. Thus this is an induced embedding of order $40$. Theorem~\ref{thm:rounded} gives $L_*(9,5)=40$ (attained at $h=5$), proving $\eta(D_{9,5})=40$.
\end{example}

\begin{theorem}[Finite arithmetic exactness range]\label{thm:finite150}
For every pair $1\le b\le a\le150$, put
\[
U(a,b)=\min_{\substack{A\ge a\\B\ge b}}\min\{3(A+B)-\gcd(A,B),\,5\max(A,B)\}.
\]
Then $\eta(D_{a,b})=L_*(a,b)=U(a,b)$. This is a finite computer-assisted theorem requiring only exact integer arithmetic, not SAT solving.
\end{theorem}
\begin{proof}
By Theorem~\ref{thm:rounded}, $L_*(a,b)\le\eta(D_{a,b})$. For $A\ge a$ and $B\ge b$, the graph $D_{a,b}$ is an induced subgraph of $D_{A,B}$; Theorem~\ref{thm:doublestar} therefore gives $\eta(D_{a,b})\le U(a,b)$. Now $U(a,b)\le f(a,b)\le750$ for $a\le150$. For $\max(A,B)\ge251$, each branch of $f(A,B)$ is at least $3\max(A,B)+2\min(A,B)\ge755$, so $U$ is obtained by enumerating $a\le A\le250$, $b\le B\le250$. The lower bound $L_*$ is a minimum over the explicitly finite range $1\le h\le a+b-1$. The complete exact-integer enumeration over the $11\,325$ admissible pairs found equality in every case; it is the script \texttt{double\_star\_bounds.py} (data availability statement), which uses only the Python standard library and runs in a few seconds. Therefore the two proved inequalities coincide for every input in the stated finite range.
\end{proof}

\begin{remark}\label{rem:oldtable}
In Table~\ref{tab:doublestars}, $\eta(D_{a,b})=f(a,b)$ except at the three entries in bold, $(9,5)$, $(11,6)$ and $(12,5)$. Table~\ref{tab:firstfailures} lists all pairs with $a+b\le20$ for which $f$ fails; the values in its $U$ column are given by Theorem~\ref{thm:finite150}.
\end{remark}

\begin{table}[ht]
\centering\small
\begin{tabular}{rrrrrr}
\toprule
$a$ & $b$ & $f(a,b)$ & $L_*(a,b)$ & $U(a,b)$ & $f-U$\\
\midrule
9&5&41&40&40&1\\
11&6&50&48&48&2\\
12&5&50&48&48&2\\
14&5&56&55&55&1\\
15&4&56&55&55&1\\
13&7&59&56&56&3\\
14&6&58&56&56&2\\
\bottomrule
\end{tabular}
\caption{All failures of the unrounded formula with $a+b\le20$, ordered by total leaf count and then by $a$. These are exact finite arithmetic consequences of Theorems~\ref{thm:rounded} and~\ref{thm:finite150}, not additional SAT claims; the script \texttt{double\_star\_bounds.py} lists them.}
\label{tab:firstfailures}
\end{table}

\begin{table}[ht]
\centering
\small
\begin{tabular}{c|cccccccc}
\toprule
$a\backslash b$ & 1 & 2 & 3 & 4 & 5 & 6 & 7 & 8\\
\midrule
1 & 5 \\
2 & 8 & 10 \\
3 & 11 & 14 & 15 \\
4 & 14 & 16 & 20 & 20 \\
5 & 17 & 20 & 23 & 25 & 25 \\
6 & 20 & 22 & 24 & 28 & 30 & 30 \\
7 & 23 & 26 & 29 & 32 & 35 & 35 & 35 \\
8 & 26 & 28 & 32 & 32 & 38 & 40 & 40 & 40 \\
9 & 29 & 32 & 33 & 38 & \textbf{40} & 42 & 45 & 45 \\
10 & 32 & 34 & 38 & 40 & 40 & 46 & 50 & 50 \\
11 & 35 & 38 & 41 & 44 & 47 & \textbf{48} & 53 & 55 \\
12 & 38 & 40 & 42 & 44 & \textbf{48} & 48 & 56 & 56 \\
\bottomrule
\end{tabular}
\caption{Values of $\eta(D_{a,b})$ for $a\le12$ and $b\le\min(a,8)$, given by Theorem~\ref{thm:finite150}; the values with $b\mid a$ also follow from Corollary~\ref{cor:dsdivides}. The entries in bold are smaller than $f(a,b)$; all others equal $f(a,b)$. The term $5a$ of $f$ is strictly smaller than $3(a+b)-\gcd(a,b)$ at $(5,4)$, $(6,5)$, $(7,6)$, $(8,7)$, $(9,7)$, $(9,8)$, $(10,8)$ and $(11,8)$, and the two terms are equal at $(a,a)$, $(4,3)$, $(7,5)$, $(8,6)$ and $(10,7)$.}
\label{tab:doublestars}
\end{table}

\begin{corollary}\label{cor:three}
For $n\ge4$, the largest value of $\eta$ over trees of diameter three on $n$ vertices is $3n-7$, attained by $D_{n-3,1}$. Consequently no constant $c<3$ satisfies $\eta(T)\le cn$ for all trees $T$.
\end{corollary}

\begin{proof}
The trees of diameter three on $n$ vertices are the double stars $D_{a,b}$ with $a+b=n-2$. By Theorem~\ref{thm:doublestar}, $\eta(D_{a,b})\le3(n-2)-\gcd(a,b)\le3n-7$, and by Corollary~\ref{cor:dsdivides}, $\eta(D_{n-3,1})=3(n-3)+2=3n-7$. Finally $(3n-7)/n\to3$.
\end{proof}

\subsection{Trees on at most eleven vertices}
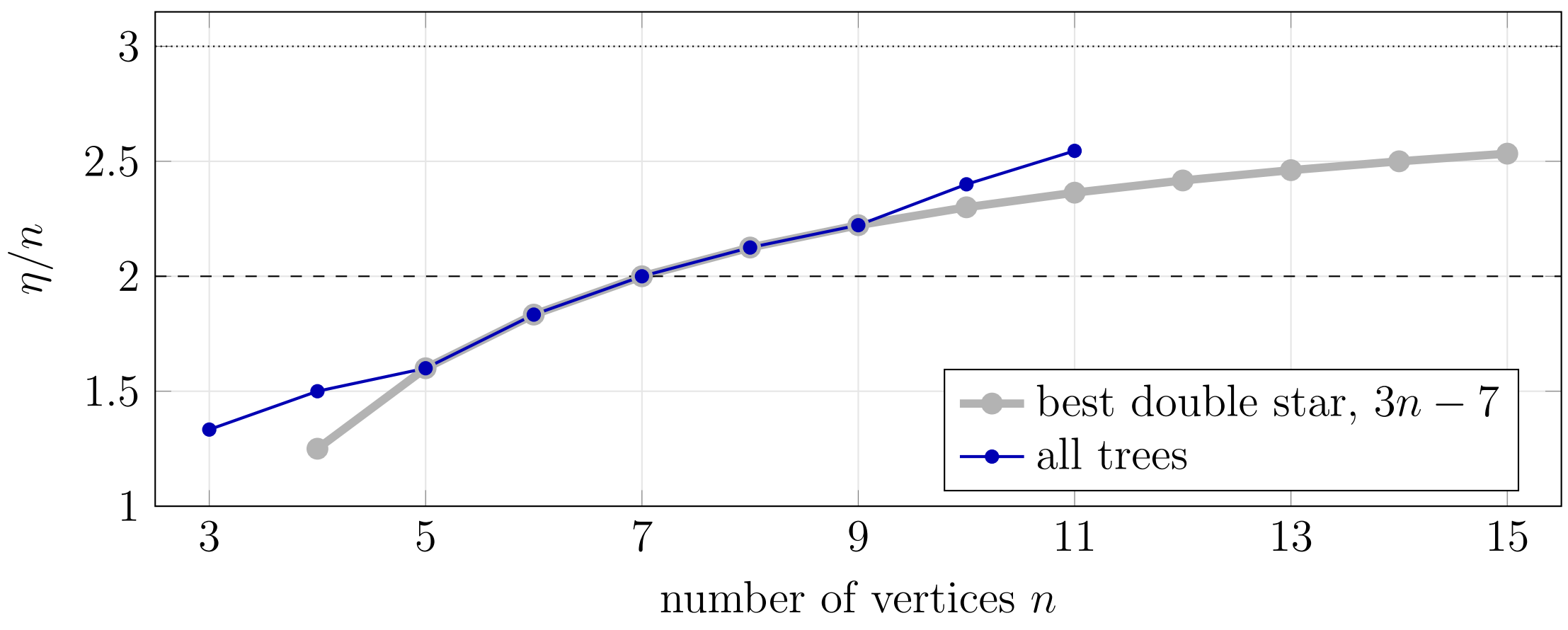
\begin{figure}[ht]
\centering
\begin{tikzpicture}
\begin{axis}[width=0.92\textwidth, height=6.2cm, xmin=2.5, xmax=15.5, ymin=1, ymax=3.15,
  xlabel={number of vertices $n$}, ylabel={$\eta/n$}, xtick={3,5,7,9,11,13,15}, ytick={1,1.5,2,2.5,3},
  legend pos=south east, legend cell align=left, grid=major, grid style={gray!20}]
\addplot[gray!60, line width=2.2pt, mark=*, mark options={fill=gray!60, draw=gray!60, scale=0.9}] coordinates {(4,1.2500) (5,1.6000) (6,1.8333) (7,2.0000) (8,2.1250) (9,2.2222) (10,2.3000) (11,2.3636) (12,2.4167) (13,2.4615) (14,2.5000) (15,2.5333)};
\addlegendentry{best double star, $3n-7$}
\addplot[blue!70!black, thick, mark=*, mark options={scale=0.75}] coordinates {(3,1.3333) (4,1.5000) (5,1.6000) (6,1.8333) (7,2.0000) (8,2.1250) (9,2.2222) (10,2.4000) (11,2.5455)};
\addlegendentry{all trees}
\addplot[black, dashed, domain=2.5:15.5] {2};
\addplot[black, densely dotted, domain=2.5:15.5] {3};
\end{axis}
\end{tikzpicture}
\caption{The largest value of $\eta(T)/n$ over all trees on $n$ vertices, computed for $n\le11$, and over double stars, where it is $(3n-7)/n$ for every $n\ge4$ by Corollary~\ref{cor:three}. The dashed line is $\eta=2n$; the dotted line at $3$ is the limit of $\eta(D_{n-3,1})/n$.}
\label{fig:trees}
\end{figure}

Table~\ref{tab:trees} and Figure~\ref{fig:trees} give, for each $n\le11$, the largest value of $\eta$ on trees with $n$ vertices and the number of trees attaining it; for $n\ge8$, the trees below the maximum are certified only to lie below it (Section~\ref{sec:computation}). For $5\le n\le 9$ a double star attains the maximum. On ten and eleven vertices, trees of diameter four do better. The unique extremal tree on eleven vertices consists of a vertex $u$ with five leaves, adjacent to a vertex $w$ that carries two pendant paths of length two; it has $\eta=28=2.55n$.

\begin{table}[ht]
\centering
\small
\begin{tabular}{rrrrrr}
\toprule
$n$ & trees & $\max\eta$ & $\max\eta/n$ & extremal trees & best double star\\
\midrule
3 & 1 & 4 & 1.33 & 1 & --\\
4 & 2 & 6 & 1.50 & 1 & 5\\
5 & 3 & 8 & 1.60 & 2 & 8\\
6 & 6 & 11 & 1.83 & 1 & 11\\
7 & 11 & 14 & 2.00 & 2 & 14\\
8 & 23 & 17 & 2.12 & 1 & 17\\
9 & 47 & 20 & 2.22 & 8 & 20\\
10 & 106 & 24 & 2.40 & 3 & 23\\
11 & 235 & 28 & 2.55 & 1 & 26\\
\bottomrule
\end{tabular}
\caption{Trees on at most eleven vertices: the largest value of $\eta$, the number of trees attaining it, and the largest $\eta(D_{a,b})$ with $a+b+2=n$ (for $n=3$ there is no double star).}
\label{tab:trees}
\end{table}

\section{A census of small graphs}\label{sec:census}

Table~\ref{tab:census} summarises the values of $\eta$ and $\ecyc$ for the 995 connected graphs on two to seven vertices of the atlas of \citet{ReadWilson1998}. Three features stand out.

\begin{itemize}
\item \emph{Extremal values.} The largest ratio $\eta/n$ grows from $1$ to $20/7\approx2.86$. The value $20$ on seven vertices is attained by two graphs whose Erd\H{o}s--Kelly floors are $11$ and $12$, respectively (atlas indices $668$ and $775$). One of them is the wheel $W_6$ (a $5$-cycle and a hub) with a pendant vertex attached to the hub: since $\eta(W_6)=10$, adding one leaf doubles $\eta$.
\item \emph{Cyclic gaps.} There are $50$ graphs with $\ecyc>\eta$, the first on five vertices (Proposition~\ref{prop:smallestgap}), and the largest ratio $\ecyc/\eta$ is $4/3$.
\item \emph{Lower bounds.} The floor $\EK$ is exact for $169$ graphs, and the linear floor $L$ for $78$ (Example~\ref{ex:blowups}).
\end{itemize}

\begin{table}[ht]
\centering
\small
\begin{tabular}{rrrrrrrr}
\toprule
$n$ & graphs & $\max\eta$ & $\max\eta/n$ & $\eta=n$ & $\ecyc>\eta$ & $\eta=\EK$ & $\eta=L$\\
\midrule
2 & 1 & 2 & 1.00 & 1 & 0 & 1 & 1\\
3 & 2 & 4 & 1.33 & 1 & 0 & 2 & 2\\
4 & 6 & 6 & 1.50 & 2 & 0 & 6 & 5\\
5 & 21 & 10 & 2.00 & 2 & 1 & 16 & 9\\
6 & 112 & 14 & 2.33 & 5 & 7 & 42 & 21\\
7 & 853 & 20 & 2.86 & 3 & 42 & 102 & 40\\
\midrule
total & 995 & & & 14 & 50 & 169 & 78\\
\bottomrule
\end{tabular}
\caption{The census of connected graphs on at most seven vertices. The column $\eta=n$ counts the abelian Cayley graphs.}
\label{tab:census}
\end{table}

The smallest graphs with a cyclic gap appear on five vertices.

\begin{proposition}\label{prop:smallestgap}
Let $K_4^+$ be $K_4$ with a pendant vertex. Then $\eta(K_4^+)=8$, with hosts $\Z_2\times\Z_4$ and $\Z_2^3$, but $\ecyc(K_4^+)=10$. No connected graph on at most four vertices has $\ecyc>\eta$, and $K_4^+$ is the only one on five vertices.
\end{proposition}

\begin{proof}
Label the pendant vertex $(0,0)$, its neighbour $(1,2)$, and the other three vertices $(0,1)$, $(0,3)$, $(1,0)$. Then $\DE=\{(0,2),(1,1),(1,2),(1,3)\}$, and the three differences between the pendant vertex and its non-neighbours, $(0,1)$, $(0,3)$ and $(1,0)$, lie outside $\DE$, together with their negatives. Lemma~\ref{lem:difference} gives $\eta\le8$, and condition (ii) of Lemma~\ref{lem:EK} gives $\eta\ge n+\Delta-\delta=8$. The values $\ecyc=10$ and the remaining statements are computational (Section~\ref{sec:computation}): $\Z_8$ and $\Z_9$ admit no induced embedding of $K_4^+$.
\end{proof}

Table~\ref{tab:named} lists some named graphs. Among them, the Petersen graph, $Q_3$, the Gr\"otzsch graph and the icosahedron all have a cyclic gap. For the Chv\'atal graph a cyclic host is optimal: the graph is the subgraph of $\Cay(\Z_{20},\{\pm2,\pm5,\pm6\})$ induced on $\{0,1,2,4,\allowbreak5,6,7,8,\allowbreak10,13,15,19\}$, and no abelian group of order at most $19$ is a host. For the regular graphs that are not Cayley graphs, condition~(v) of Lemma~\ref{lem:EK} raises the floor: to $13$ for the Petersen graph, $14$ for the icosahedron, $15$ for the Frucht graph and $16$ for the Chv\'atal graph.

\begin{table}[ht]
\centering
\small
\begin{tabular}{lrrrrlr}
\toprule
graph & $n$ & $m$ & $\EK$ & $\eta$ & optimal hosts & $\ecyc$\\
\midrule
$K_4^+$ & 5 & 7 & 8 & 8 & $\Z_2\times\Z_4$, $\Z_2^3$ & 10\\
$W_6$ with pendant at hub & 7 & 11 & 12 & 20 & $\Z_{20}$, $\Z_2\times\Z_{10}$ & 20\\
cube $Q_3$ & 8 & 12 & 8 & 8 & $\Z_2\times\Z_4$, $\Z_2^3$ & 10\\
$P_3\,\square\,P_3$ & 9 & 12 & 12 & 15 & $\Z_{15}$ & 15\\
Petersen & 10 & 15 & 10 & 16 & $\Z_4\times\Z_4$, $\Z_2^2\times\Z_4$, $\Z_2^4$ & 20\\
Gr\"otzsch & 11 & 20 & 14 & 16 & $\Z_4\times\Z_4$, $\Z_2^2\times\Z_4$, $\Z_2^4$ & 20\\
icosahedron & 12 & 30 & 12 & 20 & $\Z_2\times\Z_{10}$ & 26\\
Chv\'atal & 12 & 24 & 12 & 20 & $\Z_{20}$, $\Z_2\times\Z_{10}$ & 20\\
Frucht & 12 & 18 & 12 & 24 & $\Z_{24}$ & 24\\
$P_4\,\square\,P_4$ & 16 & 24 & 20 & 24 & $\Z_{24}$ & 24\\
\bottomrule
\end{tabular}
\caption{Some named graphs. $\EK$ is the floor of Definition~\ref{def:EK}, without condition~(v). The optimal hosts are all the abelian groups of order $\eta$ that are hosts.}
\label{tab:named}
\end{table}

\section{Computation and certification}\label{sec:computation}

\subsection*{Encoding}
For a graph $G$ and a group $\Gamma$ we use a Boolean variable $x_{v,g}$ for each $v\in V(G)$ and $g\in\Gamma$, meaning $f(v)=g$, and a variable $y_c$ for each class $c=\{g,-g\}$ of non-zero elements, meaning $c\subseteq S$. The clauses say:

\begin{itemize}
\item every vertex receives exactly one label;
\item every element is used at most once;
\item for all distinct $u,v$ and $g\ne h$, the clause $\neg x_{u,g}\vee\neg x_{v,h}\vee y_{c(g-h)}$ holds if $uv\in E(G)$, and $\neg x_{u,g}\vee\neg x_{v,h}\vee\neg y_{c(g-h)}$ holds otherwise.
\end{itemize}

The at-most-one constraints are written with the sequential counter of \citet{Sinz2005}. By Lemma~\ref{lem:difference} the formula is satisfiable exactly when $G$ has an induced embedding into some Cayley graph of $\Gamma$. It has $O(n|\Gamma|)$ variables and $O(n^2|\Gamma|^2)$ clauses.

\begin{lemma}[Symmetry breaking]\label{lem:symmetry}
Fix a vertex $v_0$ and a total order on $\Gamma$. Let $W_1,\dots,W_r$ be pairwise disjoint classes of open twins or of closed twins (vertices with equal closed neighbourhoods), none containing $v_0$. If $G$ embeds in a Cayley graph of $\Gamma$, then it has an embedding with $f(v_0)=0$ in which the labels increase within each $W_i$ along a fixed enumeration.
\end{lemma}

\begin{proof}
Translation by $-f(v_0)$ preserves $\DE$ and $\DN$. Any permutation of $V(G)$ that fixes $v_0$ and permutes each $W_i$ is an automorphism of $G$, and composing $f$ with it gives another embedding with $f(v_0)=0$. Choose the permutation that sorts the labels in each class.
\end{proof}

For each $N$ from a proven lower bound upwards, and each abelian group of order $N$ (enumerated by partitions of the prime exponents), we solve the formula with these symmetry-breaking clauses added. The least $N$ admitting a satisfying assignment is $\eta(G)$; for $\ecyc$ only cyclic groups are used. The starting point is $\EK(G)$ or, for graphs in the census, the largest certified value of $\eta$ among the connected components of the vertex-deleted subgraphs (Lemma~\ref{lem:mono}), if that is larger. The total order on $\Gamma$ is the lexicographic order of coordinate vectors in invariant-factor form, and $v_0$ is a vertex of maximum degree outside the twin classes when there is one. The solver was CaDiCaL~1.5.3 \citep{BiereEtAl2020}, called through PySAT \citep{IgnatievMorgadoMarquesSilva2018}. The solving times recorded for the instances of this paper add up to about $3.4$ minutes, and the C checker took about $9$ minutes for all their proofs.

\subsection*{Certificates}
Every upper bound is certified by an explicit labelling, which we checked with Lemma~\ref{lem:difference} by computing all $n(n-1)$ differences. For every unsatisfiable instance, CaDiCaL wrote a proof in binary DRAT format; it consists of lemmas with the reverse unit propagation (RUP) property and of deletions, so it is a DRUP proof \citep{HeuleHuntWetzler2013}. We checked these proofs with a forward RUP checker, independent of the solver, which we wrote for the purpose in two separate versions, one in C and one in Python, implementing the same algorithm with two watched literals per clause. Both versions honour deletions, never retract a literal derived at the top level (which is sound, since such literals are implied by the formula) and require the empty clause. Every proof was checked with the C version, and every proof of at most $300\,000$ lines also with the Python version; the two versions agreed on all 2437 proofs checked by both, out of 2447. In a mutation test on 20 proofs, both versions rejected every truncated proof and every proof checked against a formula from which a random half of the clauses had been removed, and they gave the same verdict on each of 100 proofs with one literal of one lemma negated (89 of these were rejected; in the others the altered lemma is still a RUP consequence). DRUP proofs can also be checked with DRAT-trim \citep{WetzlerHeuleHunt2014}. The status of each value is as follows.

\begin{itemize}
\item \emph{Census.} All $995$ graphs of Table~\ref{tab:census}: 2198 unsatisfiable instances, all with checked proofs, and 1045 satisfiable ones, each with a checked labelling. Every value of $\eta$ agrees with an earlier exhaustive search made without proof logging.
\item \emph{Named graphs.} All graphs of Table~\ref{tab:named}, starting from $\EK(G)$: 103 unsatisfiable instances, all with checked proofs. At the order $\eta(G)$ every abelian group was tested, which gives the list of optimal hosts.
\item \emph{Double stars.} The values of Table~\ref{tab:doublestars} follow from Theorem~\ref{thm:finite150}, which needs exact integer arithmetic only.
\item \emph{Trees.} Trees on at most seven vertices are part of the census. Let $M_n$ be the largest value of $\eta$ in Table~\ref{tab:trees}. For each tree $T$ on $n=8,\dots,11$ vertices we tested the orders $M_n-1,M_n-2,\dots,\EK(T)$, all abelian groups of each order, and stopped at the first host; every tree that is not extremal has a host of order less than $M_n$. For each of the $13$ extremal trees, every abelian group of order between $\EK(T)$ and $M_n-1$ was excluded by a checked proof (125 instances), and a host of order $M_n$ was found. In all, the trees gave 146 unsatisfiable instances: these 125, 16 groups excluded for trees that are not extremal before a host was found, and 5 groups of order $M_n$ that do not host the extremal tree in question. The column for double stars comes from Corollary~\ref{cor:three}.
\end{itemize}

In total $2198+103+146=2447$ formulas were unsatisfiable, with 27.7 million proof lines. The labellings, the checker outputs and the SHA-256 digests of all formulas and proofs are in the deposit described in the data availability statement. CaDiCaL is deterministic, so the deposited code regenerates the formulas and proofs, and the digests confirm that they are the ones we checked: for a sample of 40 instances from this paper and its companion \citep{FokamBitjoka2026companion}, regeneration reproduced the formula and the proof byte for byte.

\subsection*{Complexity}
Deciding $\eta(G)\le N$ is in NP: a certificate consists of the invariant factors of a group of order at most $N$ and $n$ labels, which take $O(n\log N)$ bits, and Lemma~\ref{lem:difference} checks it with $O(n^2)$ group operations. Since $\eta(G)\le p^2<4n^2$ for the least odd prime $p\ge n$ (Section~\ref{sec:constructions}), the search terminates, and each formula has $O(n^6)$ clauses.

\section{Open problems}\label{sec:problems}

\begin{problem}
Determine whether the two arithmetic functions $L_*(a,b)$ and $U(a,b)$ are equal for all positive integers $a,b$. Theorem~\ref{thm:finite150} establishes this in the finite range $1\le b\le a\le150$; no general arithmetic identity or formula for $\eta(D_{a,b})$ is claimed outside the proved cases.
\end{problem}

\begin{problem}
Determine $\eta$ for trees of diameter four, which are extremal on ten and eleven vertices.
\end{problem}

\begin{problem}
Close the gap between $\Omega(n^{6/5})$ (Corollary~\ref{cor:regular}) and $O(n^2)$ for the largest value of $\eta$ over cubic graphs on $n$ vertices.
\end{problem}

\begin{problem}
Is $\ecyc(G)/\eta(G)$ bounded? In the census the largest ratio is $4/3$.
\end{problem}

\begin{problem}
Is it NP-hard to decide whether $\eta(G)\le N$? The problem lies in NP (Section~\ref{sec:computation}). For $N=n$ it asks whether $G$ is an abelian Cayley graph. For cyclic groups this recognition problem is solvable in polynomial time \citep{EvdokimovPonomarenko2004}, so $\ecyc(G)=n$ can be decided in polynomial time; polynomial algorithms are also known for some other families of abelian groups \citep{Ryabov2021}.
\end{problem}

\begin{problem}
Is $\eta(P_r\,\square\,P_r)=r(r+2)$ for all $r\ge3$? The answer is yes for $r=3,4$ (Section~\ref{sec:constructions}), and the problem is open for every $r\ge5$. For $r\ge7$ it reduces to hosts of degree four: by condition~(i) of Lemma~\ref{lem:EK}, a host of degree $d\ge5$ has order at least $2r^2-4r(r-1)/d\ge(6r^2+4r)/5$, which exceeds $r(r+2)$ when $r>6$. For $r=5$ and $r=6$, hosts of degree five or more are not excluded by this argument.
\end{problem}

\begin{problem}
Characterise the graphs with $\eta(G)=\EK(G)$.
\end{problem}

\section*{Use of AI tools}
The authors used two AI assistants. Claude Opus 5.5 (Anthropic) ran the computations of Sections~\ref{sec:census} and~\ref{sec:computation}, prepared the deposited certificates, and helped prepare the text. ChatGPT GTP-6 Astra (OpenAI) was used in the revision to check independently Claude's computations. The authors checked all proofs and references and are responsible for the paper and its computational artifacts.

\section*{Funding}
This research received no specific grant from any funding agency in the public, commercial or not-for-profit sectors.

\section*{Competing interests}
The authors declare that they have no competing interests.

\section*{Data availability}
The integer computation behind Theorem~\ref{thm:finite150}, Tables~\ref{tab:firstfailures} and~\ref{tab:doublestars} and Example~\ref{ex:counterexample} is the script \texttt{double\_star\_bounds.py}, which uses only the Python standard library; the script and its output are ancillary files of the arXiv version of this paper (arXiv:2609.01486). For Sections~\ref{sec:census} and~\ref{sec:computation}, the deposit \citet{FokamBitjoka2026data} (DOI \DepositDOI) contains the programs that write, solve and check the formulas, every labelling, the SHA-256 digest of every formula and of every proof, and the output of the C checker on every proof and of the Python checker on every proof it checked; it also holds the same material for the companion paper \citep{FokamBitjoka2026companion}. The formulas and proofs themselves (about 0.2~GB compressed for this paper) are not in the deposit: the deposited code regenerates them, and the authors will provide them on request.

\bibliographystyle{plainnat}
\bibliography{refs}

\end{document}